\documentclass[12pt,a4paper]{article}
\usepackage{amsmath,amssymb,amsthm}
\usepackage{graphicx}
\usepackage{color}

\theoremstyle{plain}
\newtheorem{theorem}{Theorem}[section]

\newtheorem{lemma}{Lemma}[section]

\theoremstyle{remark}
\newtheorem{remark}[lemma]{Remark}

\theoremstyle{definition}

\begin{document}

\newcommand{\strain}{{\boldsymbol\varepsilon}}
\newcommand{\straind}{{\boldsymbol e}}
\newcommand{\straine}{{\boldsymbol\varepsilon^e}}
\newcommand{\plast}{{\boldsymbol\varepsilon^p}}
\newcommand{\stress}{{\boldsymbol\sigma}}
\newcommand{\inter}{{\boldsymbol\xi}}
\newcommand{\ttau}{{\boldsymbol h_d}}
\newcommand{\et}{{\boldsymbol d}^{tr}}
\newcommand{\uu}{{\boldsymbol u}}
\newcommand{\no}{\nonumber}
\newcommand{\elena}{\color{red}}
\newcommand{\debole}{\rightharpoonup}
\newcommand{\debast}{\rightharpoonup^*}
\newcommand{\real}{{\bf R}}
\newcommand{\el}{l}
\newcommand{\dive}{\hbox{div }}

                                %
\newcommand{\diffe}{D}
\title{\Large A phase transition model for the helium supercooling}

                        %

                                %
\author{\normalsize Elena Bonetti \\
{\small \it Dipartimento di Matematica} \\
{\small \it Universit\`a di Pavia, via Ferrata 1, 27100 Pavia, Italy}\\
{\scriptsize e-mail address: \tt elena.bonetti@unipv.it} \\
\normalsize  Michel Fr\'emond\\
{\small \it Dipartimento di Ingegneria Civile} \\
{\small \it Universit\`a di Roma "Tor Vergata", via del Politecnico 1, 00133 Roma, Italy}\\
{\scriptsize e-mail address: \tt michel.fremond@uniroma2.it}  }
\vspace{-1.5mm}

\date{\today}     
                                %
                              %
\maketitle
                                %
\begin{abstract}
We build a predictive theory for the evolution of mixture of helium and supercooled helium at low temperature.
The absolute temperature $\theta$ and the volume
fraction $\beta$ of helium, which is dominant at temperature larger
than the phase change temperature, are the state quantities.
The predictive theory accounts for local interactions
at the microscopic level, involving the gradient of $\beta$.
The nonlinear heat flux in the supercooled phase results
from a Norton-Hoff potential. We prove that the resulting set of
partial differential equations has solutions within a convenient analytical frame.

\end{abstract}
\vskip3mm

   \renewcommand{\theequation}{\thesection.\arabic{equation}}                             

\section
{Introduction}

We investigate a phase transition model describing the thermal behaviour of helium. The model is based on some experimental observations at low temperature between $%
\theta >0~K$ and $\theta \simeq 5~K$, as described by E. Senger in \cite{senger}. Indeed,
phase change occurs at $\theta _{c}=2.17~K$ (pressure being the atmospheric pressure): the helium, He II, is a
superfluid at temperature lower than $\theta _{c}$, the helium, He I, is a
fluid at temperature larger than $\theta _{c}$. The Gorter-Mellink model
assumes a nonlinear Fourier law based on a Norton-Hoff potential, \cite{Friaa}. In our approach this model
is upgraded by
considering microscopic motions involved in the phase change, \cite{senger}, \cite{raupp}. Thus,
the volume fraction $\beta $ of the He I is introduced, as well as the gradient of the
volume fraction, $\nabla\beta $,  to account for local
interactions, \cite{{fremondlibro}}. A description of physics together with numerical results are
given in \cite{senger}.\\
Within the small perturbation we derive a predictive theory and prove that there exist solutions of the resulting set of partial differential equations in a convenient analytical frame. The main novelty of the resulting system is represented by
the nonlinear structure of the diffusive term for the temperature (parabolic) equation. Indeed, we have the contribution of two degenerating terms, characterizing different diffusive behaviour in the two phases. Hence, the evolution equation governing the phase transition contains a multivalued operator leading internal constraints on the phase variable.

In Section 2 we detail the model derivation, leading to an initial and boundary value problem, whose analytical formulation is made precise in Section 3. Existence of a solution is proved in Section 4 by a fixed point argument, mainly combined with lower semicontinuity results.

                      %
\section
{The model}
In this paper, we aim to model the phenomenon of helium supercooling by use
of the phase transitions theory. In particular, our two phases are given by the helium at its
normal state and at the supercooled state.

We consider helium located in a smooth bounded domain
$\Omega\subset\real^n$, $n\leq3$. 
At a first step, we study the evolution of the supercooling
process during a finite time interval $(0,T)$ and assume that no macroscopic deformations
act during the phenomenon.
We use the notation $Q:=\Omega\times(0,T)$.

First, we introduce the state variables of the model $(\theta,\beta)$: $\theta$ stands
for
the absolute temperature; the phase parameter $\beta$, denoting the volume fraction of helium
at its normal state, (HE I), is
\begin{equation}\label{constr}
\beta\in[0,1].
\end{equation}
Assuming that the two phases can coexist at each point, with suitable proportions, and that
no voids nor overlapping can occur between the two phases, we let
$1-\beta$ be the volume fraction of supercooled helium, (HE II).

The free energy is stated as follows
\begin{equation}\label{freen}
\Psi(\theta,\beta,\nabla\beta)=-c_s\theta\log\theta-{\el\over{\theta_c}}(\theta-\theta_c)\beta+k|\nabla\beta|^2+I_{[0,1]}(\beta),
\end{equation}
where $\theta_c$ is the phase transition temperature between helium
and supercooled helium, $\el$ the phase change latent heat, $c_s$ the
heat capacity, and $k$ an interaction coefficient between the two phases holding at a microscopic level.
The indicator function
$I_{[0,1]}(\beta)$ accounts for \eqref{constr}, as it is
$I_{[0,1]}(\beta):=0$ if $\beta\in[0,1]$ and $I_{[0,1]}(\beta):=+\infty$,
otherwise, \cite{moreau}. Then, we introduce the dissipative variables of the model
$\beta_t$ and $\nabla\theta$,
accounting for the thermodynamical evolution of the system, and make precise the
pseudo-potential of dissipation, \cite{{fremondlibro}}
\begin{equation}\label{diss}
\Phi(\beta_t,\nabla\theta,\theta,\beta)={\mu\over2}|\beta_t|^2+\frac d\theta\left(\frac {|{\nabla\theta}|^2}2+(1-\beta)\frac{|{\nabla\theta}|^p}p\right),
    \end{equation}
where $1<p<2$ and  $\mu,d$ are positive constants. For the sake of simplicity, in the sequel we let
let $\mu=d=\el=k=c_s=1$.

The balance equations are recovered from the classical laws of
continuum thermo-mechanics. More precisely, we exploit a generalized version
of  the principle of virtual
power including also the effects of microscopic motions which are
responsible for the phase transition, \cite{{fremondlibro}}. Thus, we can recover the equation
governing the evolution
of the phase parameter $\beta$ as a balance
equation for microscopic movements. It is (${\bf
  n}$ is the normal vector to the boundary)
\begin{equation}\label{eqbilan}
B-\dive{\bf H}=0\hbox{ in }\Omega,\quad{\bf H}\cdot{\bf n}=0\hbox{ on }\partial\Omega,
\end{equation}
where $B$ and ${\bf H}$ are new interior forces for which we are
specifying in a moment the  constitutive relations. Note that we are not considering applied
volume or surface forces acting on the microscopic level.
The second equation is given by the energy balance, mainly governing
the evolution of the temperature
\begin{equation}\label{eqener}
e_t+\dive{\bf q}=r+B\beta_t+{\bf H}\cdot\nabla\beta_t\hbox{ in }
\Omega,\quad{\bf q}\cdot{\bf n}=0\hbox{ on }\partial\Omega,
\end{equation}
$e$ denoting the internal energy, ${\bf q}$ the heat flux, $r$ an exterior heat source, and
$B\beta_t+{\bf H}\cdot\nabla\beta_t$ microscopic mechanically induced heat sources.

Now, let us specify the constitutive relations for the above involved
physical quantities. As usual, the internal energy  $e$ is specified  by
\begin{equation}\label{ener}
e=\Psi+\theta s,
\end{equation}
where  the
entropy $s$ is
\begin{equation}\label{entropia}
s=-{{\partial\Psi}\over{\partial\theta}}.
\end{equation}
Then, $B$ is given by the sum of a non-dissipative contribution (derived by the
free energy $\Psi$) and a
dissipative one (related to the pseudo-potential $\Phi$), i.e.
\begin{equation}\label{Bi}
B=B^{nd}+B^d={{\partial\Psi}\over{\partial\beta}}+{{\partial\Phi}\over{\partial\beta_t}},
\end{equation}
while  ${\bf H}={\bf H}^{nd}$ is  taken only as a non-dissipative vector defined by
\begin{equation}\label{acca}
{\bf H}^{nd}={{\partial\Psi}\over{\nabla\beta}}.
\end{equation}
Finally, the heat flux is recovered by the pseudo-potential of
dissipation through the following relation
\begin{equation}\label{flusso}
{\bf q}=-\theta{{\partial\Phi}\over{\nabla\theta}}.
\end{equation}
Let us point out that, due to the explicit form of $\Phi$, we eventually get
\begin{equation}\label{espliflu}
{\bf q}=-\beta\nabla\theta-(1-\beta)|{\nabla\theta}|^{p-2}\nabla\theta.
\end{equation}
In particular, we find the classical Fourier heat flux law for the normal state helium ($\beta=1$),
while the supercooled helium ($\beta=0$) is characterized by a lower order diffusion
term for the temperature ($p<2$).
%
%

To derive our PDE system,
we substitute \eqref{entropia}-\eqref{flusso} (cf. \eqref{espliflu}) in \eqref{eqbilan}-\eqref{eqener}. On
account of \eqref{freen} and \eqref{diss}, assuming the small perturbations assumption
(so that some higher order dissipative nonlinearities are neglected) and letting the model derived close
to the phase transition temperature $\theta_c$, we recover the following equations, with $1<p<2$, (cf. Remark \ref{approssimo1})
\begin{align}\label{eqIo}
&\theta_t+\beta_t-\dive(\beta\nabla\theta+(1-\beta)|{\nabla\theta}|^{p-2}\nabla\theta)=r\\\label{eqIIo}
&\beta_t-\Delta\beta+\partial
I_{[0,1]}(\beta)\ni{1\over{\theta_c}}(\theta-\theta_c).
\end{align}
Then, we have to combine \eqref{eqIIo}
with suitable initial and boundary conditions. In particular, Cauchy conditions hold in $\Omega$
\begin{equation}\label{iniz}
\theta(0)=\theta_0,\quad\beta(0)=\beta_0,
\end{equation}
and on $\partial\Omega$, $\partial_n$ being the normal derivative operator (cf. \eqref{eqener} and \eqref{eqbilan})
\begin{equation}\label{bound}
\partial_n\chi=0,\quad(\beta\nabla\theta+(1-\beta)|{\nabla\theta}|^{p-2}\nabla\theta)\cdot{\bf n}=0
\end{equation}
\begin{remark}\label{approssimo1}
Note that actually \eqref{eqIo} is obtained
regularizing the right hand side of the following equation
\begin{equation}
\label{eqIino}
\theta_t+{\theta\over{\theta_c}}\beta_t-\dive(\beta\nabla\theta+(1-\beta)|{\nabla\theta}|^{p-2}\nabla\theta)
=|{\beta_t}|^2+r \end{equation}
neglecting the dissipative higher
order nonlinearity (within the small perturbations assumption) and letting $\theta/\theta_c\sim1$
\end{remark}

\section
{Analytical formulation}

In this section, we make precise the abstract version of the analytical problem we are dealing with \eqref{eqIo}-\eqref{bound},
and state the main existence result (obtained under suitable assumptions on the data).

Concerning the regularity of initial data, we let
\begin{align}\label{regthetao}
&\theta_0\in L^2(\Omega), \\\label{regbetao}
&\beta_0\in H^1(\Omega),\quad \beta_0\in[0,1]\hbox{ a.e. in
}\Omega.
\end{align}
Hence, we take in \eqref{eqIo}
\begin{equation}\label{regr}
r\in L^2(0,T;L^2(\Omega)).
\end{equation}
Now, we can state the main existence result, which actually refers to an abstract version of the system
\eqref{eqIIo}, \eqref{eqIo} combined with \eqref{bound}, \eqref{iniz}.

We first clarify the abstract  setting we need as well as some related notation. The Hilbert triplet
\begin{equation}\label{triplet}
V:=H^1(\Omega)\hookrightarrow  H:=L^2(\Omega)\hookrightarrow  V',
\end{equation}
is introduced
$H$ being identified with its dual space, as usual. Then, by $\|{\cdot}\|_X$ we denote the norm in a Banach space
$X$ and in any power of it. By $\langle\cdot,\cdot\rangle$ we denote the duality pairing between $V',V$.

Then, we make precise the notion of solution for the above problem. We are looking for
\begin{align}\label{regbasetheta}
&\theta\in H^1(0,T;V')\cap L^2(0,T;H)\cap L^p(0,T;W^{1,p}(\Omega)),\\\label{regbeta}
&\beta\in H^1(0,T;H)\cap L^\infty(0,T;V)\cap L^2(0,T;H^2(\Omega)),\\\label{regxi}
&\xi\in L^2(0,T;H),
\end{align}
with
\begin{align}\label{xidiff}
&\xi\in \partial I_{[0,1]}(\beta)\quad\hbox{a.e. in }Q,\\\label{betatheta}
&\beta^{1/2}\nabla\theta\in L^2(Q),
\end{align}
fulfilling a.e. in $(0,T)$ in $V'$

\begin{align}\label{eqI}
&\theta_t+\beta_t-\dive(\beta\nabla\theta+(1-\beta)|{\nabla\theta}|^{p-2}\nabla\theta)
=r\\\label{eqII}
&\beta_t-\Delta\beta+\xi={1\over{\theta_c}}(\theta-\theta_c).
\end{align}
Here $-\hbox{div }(H)^3\rightarrow V'$ stands for the abstract operator
$$
\langle-\hbox{div {\bf v}},\phi\rangle=\int_\Omega {\bf v}\cdot\nabla\phi,
$$
and $-\Delta: V\rightarrow V'$ for
$$
\langle-\Delta v,\phi\rangle=\int_\Omega \nabla v\cdot\nabla\phi.
$$

The following theorem holds.

\begin{theorem}\label{esistenza}
Let \eqref{regthetao}-\eqref{regr} hold. Then, there exists a solution $(\theta,\beta,\xi)$  to \eqref{eqI}, \eqref{eqII}, \eqref{xidiff}, \eqref{iniz} with regularity \eqref{regbasetheta}-\eqref{regxi}, \eqref{betatheta}. In particular, \eqref{eqII} is actually solved a.e. in $Q$.
\end{theorem}

%
%

\section
{Existence result}

In order to prove the existence result stated  by Theorem \ref{esistenza} we apply a fixed point argument
on a regularized version of our problem, showing that it admits, at least locally in time,
existence of a solution.
Hence, exploiting an {\sl a priori estimates-passage to the limit} procedure we pass
to the limit showing that the limit problem admits a solution, at least in a weak sense. We point out that, actually,
the limit problem is solved globally in time, due to suitable a priori estimates on the solutions,
allowing to extend them on the whole time interval.

\subsection{The regularized problem}

First, letting $\varepsilon>0$, we regularize \eqref{eqI} as follows
\begin{equation}\label{eqIe}
\theta_t+\beta_t-\varepsilon\Delta\theta-\dive(\beta\nabla\theta+(1-\beta)|{\nabla\theta}|^{p-2}\nabla\theta)=r,
\end{equation}
and combine it with \eqref{eqII} and \eqref{iniz}.
Then,
we construct an operator ${\cal T}$ fulfilling the assumption of the Schauder theorem and such that its fixed points are eventually solutions to our system \eqref{eqIe},\eqref{eqII}, \eqref{xidiff} with \eqref{iniz} .
To this aim, we introduce
\begin{equation}
{\cal B}:=\{v\in L^2(0,\widehat T;H):\,\|{v}\|_{L^2(0,\widehat T;H)}\leq R\},
\end{equation}
here $R>0$ is fixed and $\widehat T\in(0,T]$ will be chosen later.
Note that, from now on, we denote by $c$ possibly different positive constants not depending on the
solutions, but just on the data of the problem, and by $c(R)$ positive constants depending in particular on $R$.

First, we take $\bar\theta\in{\cal B}$ and substitute $\theta$ in \eqref{eqII}.
Fairly standard results in the theory of evolution equations associated with maximal monotone
operators ensure the existence and uniqueness of
$$
\beta:={\cal T}_1(\bar\theta)\in H^1(0,\widehat T;H)\cap L^2(0,\widehat T;H^2(\Omega))
$$
solving the resulting equation (actually a.e.), with \eqref{iniz} and \eqref{bound} (note that the right
hand side of \eqref {eqII} belongs to
$L^2(\Omega\times(0,\widehat T))$ and the regularity of the initial datum is given by \eqref{regbetao}).
Then, we can exploit the following a priori estimate on the solution $\beta$.
We test \eqref{eqII} by $\beta_t-\Delta\beta$
and integrate over $(0,t)$, $t\in(0,\widehat T)$.
Integrating by parts, using a generalization of the chain rule for subdifferential
operators and exploiting the monotonicity of the
subdifferential $\partial I_{[0,1]}(\beta)$ (cf. \cite{brezis} and \cite{bonetti1}), as well aas  the positivity
of the indicator function $I_{[0,1]}(\cdot)$, we get (cf. also \eqref{regbetao})
\begin{align}\label{stI}
&\|{\beta_t}\|^2_{L^2(0,t;H)}+\|{\nabla\beta(t)}\|^2_H+\|{\Delta\beta}\|^2_{L^2(0,t;H)}+\int_\Omega I(\beta(t))\\\no
&\leq \|{\beta_0}\|^2_V+c\int_0^t\|{\bar\theta-\theta_c}\|_H\left(\|{\beta_t}\|_H+\|{\Delta\beta}\|_H\right)\\\no
&\leq c\left(1+\|{\bar\theta}\|^2_{L^2(0,\widehat T;H)}\right)+{1\over2}\|{\beta_t}\|^2_{L^2(0,t;H)}+
{1\over2}\|{\Delta\beta}\|^2_{L^2(0,t;H)},
\end{align}
where we have used the Young inequality.
By \eqref{stI}, using the Gronwall lemma,  we eventually deduce
\begin{equation}\label{boundI}
\|{\beta}\|_{H^1(0,\widehat T;H)\cap L^\infty(0,\widehat T;V)\cap L^2(0,\widehat T;H^2(\Omega))}\leq c(R).
\end{equation}
Moreover, we can infer that $\beta\in[0,1]$ a.e.. By a comparison in \eqref{eqII}  $\xi\in \partial I(\beta)$
satisfies
\begin{equation}\label{boundIxi}
\|\xi\|_{L^2(0,\widehat T;H)}\leq c(R).
\end{equation}

Now, as a second step, we take $\beta={\cal T}_1(\bar\theta)$ in \eqref{eqIe} and look for a corresponding solution
$\theta={\cal T}_2(\beta)$ (actually depending on $\varepsilon>0$).
To this aim, we apply the theory of parabolic nonlinear evolution equations (recall that $\beta\in[0,1]$), and get existence
and uniqueness of a solution (see \cite{lions})
$$
\theta\in H^1(0,\widehat T;V')\cap L^2(0,\widehat T;V).
$$
Then, let us proceed detailing some a priori estimates on this solution. For the sake of simplicity,
if it is allowed, we directly perform estimates in which the positive bound $c$ does not depend on the approximating parameter $\varepsilon$, as we are interested, in a second step, to pass to the limit on $\varepsilon\searrow0$ on the whole $(0,T)$.
We test  \eqref{eqIe} by $\theta$ and integrate over $(0,t)$, $t\in[0,T]$.
After integrating by parts in time and exploiting the H\"older
inequality, we get
\begin{align}\label{stimaI}
&{1\over2}\|{\theta(t)}\|^2_{H}-{1\over2}\|{\theta_0}\|_{H}+\varepsilon\int_0^t\|\nabla\theta\|^2_H
+\int_0^t\int_\Omega\beta|{\nabla\theta}|^2+\int_0^t\int_\Omega|{\nabla\theta}|^p\\\no
&\leq \int_0^t\|{r}\|_H\|{\theta}\|_H+\int_0^t\|\beta_t\|_H\|\theta\|_H+\int_0^t\int_\Omega\beta|{\nabla\theta}|^p.
\end{align}
Now, let us handle the last term on the right hand side of \eqref{stimaI} exploiting the Young inequality
\begin{align}\label{stimaIbis}
&\int_0^t\int_\Omega\beta|{\nabla\theta}|^p=\int_0^t\int_\Omega\beta^{1-p/2}\beta^{p/2}|{\nabla\theta}|^p\\\no
&\leq
\delta\int_0^t\int_\Omega\left(\beta^{p/2}|{\nabla\theta}|^p\right)^{2/p}+C_\delta\int_0^t\int_\Omega|{\beta^{1-p/2}}|^{\frac 2{2-p}}\\\no
&\leq \delta\int_0^t\int_\Omega\left(\beta|{\nabla\theta}|^2\right)+c\left(1+\|{\beta}\|^2_{L^2(0,\widehat T;H)}\right).
\end{align}
Thus, letting, e.g., $\delta=1/2$, we combine \eqref{stimaIbis} with
\eqref{stimaI}, and  the Gronwall lemma implies
\begin{equation}\label{boundII}
\|{\theta}\|_{L^\infty(0,\widehat T;H)}+\|{\beta^{1/2}\nabla\theta}\|_{L^2(\Omega\times(0,\widehat T))}
+\|{\nabla\theta}\|_{L^p(\Omega\times(0,\widehat T))}+\varepsilon^{1/2}\|\nabla\theta\|_{L^2(\Omega\times(0,\widehat T))}
\leq c(R).
\end{equation}
Then, a comparison in \eqref{eqI} yields (at least)
\begin{equation}\label{boundIII}
\|{\theta}\|_{H^1(0,\widehat T;V')}\leq c.
\end{equation}

Now, we define
$$
{\cal T}(\bar\theta)={\cal T}_2({\cal T}_1(\bar\theta)).
$$
It results that ${\cal T}:{\cal B}\rightarrow{\cal B}$ is well-defined,
at least for some small $\widehat T$. Indeed, by virtue
of \eqref{boundII} we can infer that
$$
\|\theta\|_{L^2(0,\widehat T;H)}=\left(\int_0^{\widehat T}\|\theta\|^2_H\right)^{1/2}
\leq \|\theta\|_{L^\infty(0,\widehat T;H)}\frac {\widehat T^{1/2}} 2\leq \frac {c(R)}2\widehat T^{1/2}
$$
and the right hand side is less than $R$ if $\widehat T$ is sufficiently small.
However, as the above estimates
do not depend on $\widehat T$. Then, as  $\theta\in C^0([0,\widehat T];H)$, $\beta\in C^0([0,\widehat T];V)$,
due to \eqref{regthetao}-\eqref{regbetao}, we will be able to extend the result on the whole time interval $(0,T)$.
Thus, for the sake of simplicity, in the sequel  we directly refer to the interval $(0,T)$.
Hence, it is clear that any fixed point of ${\cal T}$  is a solution  $\theta,{\cal T}_1(\theta))$ to \eqref{eqIe},\eqref{eqII}, \eqref{iniz}.

Now, to prove that ${\cal T}$ admits a fixed point exploiting the Schauder theorem, we have to show that it is compact and continuous w.r.t.
to the topology of $L^2(0,T;H)$. Compactness easily follows by \eqref{boundII} and \eqref{boundIII}
(see \cite{simon}), as $\varepsilon>0$ (we have that $L^2(0,T;H)$ is compact in $H^1(0,T;V')\cap L^2(0,T;V)$.

Now, we aim to prove (strong) continuity.
To this aim we take a sequence
\begin{equation}\label{contI}
\bar\theta_n\rightarrow\bar\theta\quad\hbox{in }L^2(0,T;H)
\end{equation}
and show that the corresponding $\theta_n$ strongly converges to
\begin{equation}
\theta_n={\cal T}(\bar\theta_n)\rightarrow\theta={\cal T}(\bar\theta)
\end{equation}
strongly in $L^2(0,T;H)$.

First, we observe that \eqref{boundI}, \eqref{boundII}, and \eqref{boundIII} hold for $\beta_n={\cal T}_1(\bar\theta_n)$ and
$\theta_n$ for constants $c$ independent of $n$. Thus, by weak and weak star compactness results, we get,
at least for some suitable subsequences,
\begin{equation}\label{convI}
\beta_n\debast\beta\quad\hbox{in }H^1(0,T;H)\cap L^\infty(0,T;V)\cap L^2(0,T;H^2(\Omega)).
\end{equation}
Hence, we use strong compactness theorems yielding (at least for subsequences)
\begin{equation}\label{convII}
\beta_n\rightarrow\beta\quad\hbox{in }C^0([0,T];H)\cap L^2(0,T;H^{2-\delta}(\Omega)),\quad\delta>0.
\end{equation}
As $\beta_n\in[0,1]$ a.e. and, for a subsequence, $\beta_n\rightarrow\beta$ a.e. (cf. \eqref{convII}), the Lebesgue theorem
ensures that
\begin{equation}\label{convIII}
\beta_n\rightarrow\beta\quad\hbox{in }L^q(Q),\quad\forall q<+\infty.
\end{equation}
Analogously, due to \eqref{boundIxi}, we can infer that
\begin{equation}\label{convIV}
\xi_n\rightharpoonup \xi\quad\hbox{in }L^2(0,T;H),
\end{equation}
so that \eqref{convII} and \eqref{convIV} lead to the identification (cf. \cite{brezis})
\begin{equation}\label{convV}
\xi\in\partial I(\beta),
\end{equation}
a.e. in $Q$.
The above convergences \eqref{convI}-\eqref{convV}  allow us to pass to the limit in the equation \eqref{eqII}
as $n\rightarrow+\infty$, actually in $L^2(0,T;H)$ and thus to identify the limit equation a.e. in $Q$.
By the uniqueness of the limit equation, once $\bar\theta$ is fixed, we identify  $\beta={\cal T}_1(\bar\theta)$
and extend \eqref{convI}-\eqref{convV} to the whole sequences.

Now, we deal with \eqref{eqI} where $\beta_n$ is fixed and look for suitable convergence results of the corresponding
solutions  $\theta_n={\cal T}_2(\beta_n)$, as $n$ tends to $+\infty$.
We first point out that \eqref{boundII} and \eqref{boundIII}, due to weak and weak star compactness results,
imply that the following
convergence holds, up to the extraction of some suitable subsequences,
\begin{align}\label{convVI}
&\theta_n\debast\theta\quad\hbox{in }{H^1(0,T;V')\cap
L^\infty(0,T;H)\cap L^p(0,T;W^{1,p}(\Omega))},\\\label{convVII}
&\beta^{1/2}_n\nabla\theta_n\debole\eta\quad\hbox{in }L^2(0,T;H),\\\label{convVIIbis}
&\varepsilon^{1/2}\theta_n\debole\varepsilon^{1/2}\theta\quad\hbox{in }L^2(0,T;V).
\end{align}
As $\varepsilon>0$, \eqref{boundII} implies that, at least for a subsequence,
\begin{equation}\label{convfortete}
\theta_n\rightarrow\theta\quad\hbox{in }L^2(0,T;H).
\end{equation}
Then,
we can easily identify $\eta=\beta^{1/2}\nabla\theta$ in \eqref{convVII} due to \eqref{convIII} and \eqref{convVI}. Indeed, $\beta_n$ strongly converges to $\beta$, and thus it converges almost everywhere. As a consequences $\beta_n^{1/2}\rightarrow\beta^{1/2}$ a.e. and  $\beta_n^{1/2}\in[0,1]$ a.e. (i.e. they are uniformly bounded). The Lebesgue theorem ensures that $\beta_n^{1/2}\rightarrow\beta^{1/2}$  strongly in $L^r(Q)$ for any $1\leq r<+\infty$. Finally, note that, due to \eqref{convVIIbis}, $\nabla\theta_n$ converges weakly in $L^2(Q)$.

Then, let us point out that \eqref{boundII} yields
\begin{equation}\label{boundIV}
\|{|{\nabla\theta_n}^{p-2}\nabla\theta_n|}\|_{L^{p'}(Q)}\leq
c,\quad{1\over p}+{1\over{p'}}=1,
\end{equation}
so that  \eqref{convIII}, combined with \eqref{boundIV}, yields
\begin{equation}\label{convVIII}
(1-\beta_n)|{\nabla\theta_n}^{p-2}|\nabla\theta_n\debole(1-\beta)\psi\quad\hbox{
in }L^{p'}(Q),
\end{equation}
(and consequently in  $L^2(Q)$) where
\begin{equation}
|{\nabla\theta_n}^{p-2}|\nabla\theta_n\debole\psi\quad\hbox{in }L^{p'}(Q).
\end{equation}

By the previous convergences we are allowed to pass to the limit as $n\rightarrow+\infty$
in \eqref{eqIe} and get
\begin{equation}\label{finastr}
\theta_t-\varepsilon\Delta\theta-\dive(\beta\nabla\theta+(1-\beta)\psi)+\beta_t=r,
\end{equation}
in $V'$, for a.e. $t$.
Now, our aim is to
identify $\psi$ in \eqref{finastr}  with $|{\nabla\theta}^{p-2}|\nabla\theta$.

Once $\beta$ is fixed we introduce the function
\begin{align}
&\widehat J_\beta:L^p(\Omega)\rightarrow(-\infty,+\infty],\,\widehat J_\beta(\nabla\theta)=\int_\Omega \frac{(1-\beta)}p|\nabla\theta|^p\\\no
&\hbox{ if }|\nabla\theta|^p\in L^1(\Omega),\,\widehat J_\beta(\nabla\theta)=+\infty\hbox{ otherwise}.
\end{align}
Then, the subdifferential of $\widehat J_\beta$ can be standardly computed as
\begin{equation}
\partial\widehat J_\beta=(1-\beta)|\nabla\theta|^{p-2}\nabla\theta.
\end{equation}
Our goal is to prove that $(1-\beta)\psi\in\partial\widehat J_\beta(\nabla\theta)$ a.e. in $Q$, i.e. (by definition of the subdifferential)
\begin{equation}
\int_Q(1-\beta)\psi(\nabla w-\nabla\theta)\leq\int_Q\frac{(1-\beta)}p|\nabla w|^p-\int_Q\frac{(1-\beta)}p|\nabla\theta|^p,
\end{equation}
for any $w\in L^p(0,T;W^{1,p}(\Omega))$.
To this aim we first prove that
\begin{equation}
\limsup_{n\rightarrow+\infty}\int_Q(1-\beta_n)|\nabla\theta_n|^{p-2}\nabla\theta_n\leq \int_Q(1-\beta)\psi\nabla\theta.
\end{equation}
Indeed, testing \eqref{eqI} written for $n$ by $\theta_n$ and integrating over $(0,T)$ leads to
\begin{align}\label{sci1}
&\limsup_{n\rightarrow+\infty}\int_Q(1-\beta_n)|\nabla\theta_n|^{p-2}\nabla\theta_n=\limsup_{n\rightarrow+\infty}-\frac 12
\int_\Omega|\theta_{n}(t)|^2+\frac 1 2\int_\Omega|\theta_0|^2\\\no
&-\varepsilon\int_Q|\nabla\theta_n|^2-\int_Q\beta_n|\nabla\theta_n|^2-\int_Q\theta_n\beta_{nt}+\int_Qr\theta_n\\\no
&\leq-\frac 1 2\int_\Omega|\theta(t)|^2+\frac 12\int_\Omega|\theta_0|^2-\varepsilon\int_Q|\nabla\theta|^2-\int_Q\beta|\nabla\theta|^2-\int_Q\theta\beta_{t}+\int_Qr\theta\\\no
&=+\int_Q(1-\beta)\psi\nabla\theta.
\end{align}
To prove \eqref{sci1} we have used (weak) lower semicontinuity of norms, convergences \eqref{convVI}-\eqref{convVIIbis},
\eqref{convVIII}, and \eqref{convfortete}. Hence the right-hand side is identified
by virtue of  \eqref{finastr}. As a consequence, we immediately have that
\begin{equation}\label{sci2}
\int_Q(1-\beta)\psi(\nabla w-\nabla\theta)\leq\liminf_{n\rightarrow+\infty}\int_Q(1-\beta_n)|\nabla\theta_n|^{p-2}\nabla\theta_n(\nabla w-\nabla\theta_n),
\end{equation}
for any $w\in L^p(0,T;W^{1,p}(\Omega))$. Thus, it is now a standard matter to prove the following chain of inequalities
\begin{align}\label{sci3}
&\int_Q(1-\beta)\psi(\nabla w-\nabla\theta)\leq\liminf_{n\rightarrow+\infty}\int_Q(1-\beta_n)|\nabla\theta_n|^{p-2}\nabla\theta_n(\nabla w-\nabla\theta_n)\\\no
&\leq\liminf_{n\rightarrow+\infty}\int_Q\frac{(1-\beta_n)}p(|\nabla w|^p-|\nabla\theta_n|^p)\\\no
&\leq\limsup_{n\rightarrow+\infty}\int_Q\frac{(1-\beta_n)}p(|\nabla w|^p-|\nabla\theta_n|^p)\\\no
&\leq\int_Q(1-\beta)\frac{|\nabla w|^p}p-\liminf_{n\rightarrow+\infty}\int_Q(1-\beta_n)\frac{|\nabla\theta_n|^p}p\\\no
&\leq \int_Q(1-\beta)\frac{|\nabla w|^p}p-\int_Q(1-\beta)\frac{|\nabla\theta|^p}p
\end{align}
which concludes our proof. Indeed, note in particular that, as \eqref{convVI} and \eqref{convIII} hold, we can infer that
\begin{equation}
 (1-\beta_n)^{1/p}\nabla\theta_n\debole(1-\beta)^{1/p}\nabla\theta\quad\hbox{ in }L^p(Q),
\end{equation}
so that, by weak lower semicontinuity of norms,
$$
\liminf_{n\rightarrow+\infty}\int_Q(1-\beta_n)\frac{|\nabla\theta_n|^p}p\geq\int_Q(1-\beta)\frac{|\nabla\theta|^p}p.
$$

\subsection{The existence result: the limit as $\varepsilon\searrow0$}
\newcommand\betae{\beta_\varepsilon}
\newcommand\thetae{\theta_\varepsilon}
\newcommand\xie{\xi_\varepsilon}

We are now in the position of proving Theorem \ref{esistenza} by passing to the limit in \eqref{eqIe} as
$\varepsilon\searrow0$. To this aim, after denoting by $(\thetae,\betae,\xie)$ a solution to the
system \eqref{eqIe}, \eqref{eqII}, \eqref{xidiff}, with $\varepsilon>0$, we perform the analogous
estimate as \eqref{stI} and \eqref{stimaI}-\eqref{stimaIbis}, i.e. we test \eqref{eqII} by
$\partial_t\beta_\varepsilon-\Delta\beta_\varepsilon$ and \eqref{eqIe} by $\thetae$.
After integrating over $(0,t)$ we get
\begin{align}
 &\|\partial_t\betae|^2_{L^2(0,t;H)}+\|\nabla\betae(t)\|^2_H+\|\Delta\betae|^2_{L^2(0,t;H)}+\int_\Omega I(\betae(t))\\\no
 &+\frac 1 2\|\thetae(t)\|^2_H+\varepsilon\|\nabla\thetae\|^2_{L^2(0,t;H)}+\int\int_Q\betae|\nabla\thetae|^2
 +\int\int_Q|\nabla\thetae|^2\\\no
 &\leq \frac 1 2\|\partial_t\betae|^2_{L^2(0,t;H)}+\frac 12\|\Delta\thetae\|^2_{L^2(0,t;H)}\\\no
 &+c\left(1+\|\thetae\|^2_{L^2(0,t;H)}+\int_0^t\|r\|_H\|\theta\|_H\right)+\int\int_Q\beta|\nabla\thetae|^p\\\no
 &\leq \frac 1 2\|\partial_t\betae|^2_{L^2(0,t;H)}+\frac 12\|\Delta\thetae\|^2_{L^2(0,t;H)}
 +\frac 1 2\int\int_Q\betae|\nabla\thetae|^2\\\no
 &+c\left(1+\|\thetae\|^2_{L^2(0,t;H)}+\int_0^t\|r\|_H\|\theta\|_H+\|\betae\|^2_{L^2(0,t;H)}\right).
\end{align}

Then by using the Gronwall lemma,
we get the analogous of \eqref{boundI}, \eqref{boundIxi}, \eqref{boundII}, \eqref{boundIII}, now independently of $\varepsilon$,
i.e.
\begin{align}
&\|\betae\|_{H^1(0,T;H)\cap L^\infty(0,T;V)\cap L^2(0,T;H^2(\Omega))}\leq c\\
&\|\thetae\|_{H^1(0,T;V')\cap L^\infty(0,T;H)\cap L^p(0,T;W^{1,p}(\Omega))}\leq c\\
&\|\xie\|_{L^2(0,T;H)}\leq c\\
&\|\betae^{1/2}\nabla\thetae\|_{L^2(0,T;H)}\leq c\\
&\varepsilon^{1/2}\|\thetae\|_{L^2(0,T;H)}\leq c.
\end{align}
As in the previous section (cf. \eqref{convI}-\eqref{convIV}), by compactness, we can deduce that, at least for some subsequences,
\begin{align}\label{convIe}
&\betae\debast\beta\quad\hbox{in }H^1(0,T;H)\cap L^\infty(0,T;V)\cap L^2(0,T;H^2(\Omega))\\\label{convIIe}
&\betae\rightarrow\beta\quad\hbox{in }C^0([0,T];H)\cap L^2(0,T;H^{2-\delta}(\Omega)),\quad\delta>0\\\label{convIIIe}
&\beta_n\rightarrow\beta\quad\hbox{in }L^q(Q),\quad\forall q<+\infty.\\
\label{convIVe}&\xi_n\rightharpoonup \xi\quad\hbox{in }L^2(0,T;H)
\end{align}
and (at least)
\begin{align}\label{deboletfin}
&\thetae\debast\theta\quad\hbox{in }H^1(0,T;V')\cap L^\infty(0,T;H)\cap L^p(0,T;W^{1,p}(\Omega))\\\label{fortetfin}
&\thetae\rightarrow\theta\quad\hbox{in }L^\infty(0,T;L^q(\Omega)),
\end{align}
where $q$ is such that $L^q$ is compact embedded in $W^{1,p}$, i.e. $p>\frac{nq}{q+n}$
($n$ standing for the dimension of $\Omega$). Note in particular that \eqref{convfortete} still holds
for $n=2$ for any $p\in(1,2)$,
while for $n=3$ it is ensured just for $p>6/5$.
Then, by virtue of \eqref{convIe}-\eqref{convIVe} and \eqref{deboletfin} we can pass to the limit in
\eqref{eqII} (written for $\varepsilon>0$) as $\varepsilon$ tends to zero,
identifying $\xi\in\partial I_{0,1}(\beta)$.
Hence, we are interested in pass to the limit in \eqref{eqIe}. We make use of convergences \eqref{convIIe}, \eqref{convIe},
\eqref{deboletfin}, observing in particular that
$$
-\varepsilon\Delta\thetae\rightarrow0\quad\hbox{ in }L^2(0,T;V').
$$
Arguing as in the previous section, we can infer that
\begin{equation}
\betae^{1/2}\nabla\thetae\rightharpoonup\beta^{1/2}\nabla\theta\quad\hbox{in }L^2(0,T;H),
\end{equation}
and
\begin{equation}
(1-\betae)|\nabla\thetae|^{p-2}\nabla\thetae\rightharpoonup(1-\beta)\psi\quad\hbox{in }L^{p'}(Q).
\end{equation}
We aim to identify $\psi=|\nabla\theta|^{p-2}\nabla\theta$. To this aim we proceed as in the previous section,
thus we do not detail inequalities. We just make some comment on the equivalent version of \eqref{sci1}.
Indeed, in the general case we cannot infer that
$\int_Q\thetae\partial_t\betae\rightarrow\int_Q\theta\partial_t\beta$. Thus, to prove that
$$
\limsup_{\varepsilon\searrow0}-\int_Q\thetae\partial_t\betae\leq-\int_Q\theta\partial_t\beta
$$
we use semicontinuity for the equation \eqref{eqII} (written for $\varepsilon>0$) formally tested by $\theta_c\partial_t\beta_\varepsilon$.
There holds
\begin{align}
&\limsup_{\varepsilon\searrow0}-\int_Q\thetae\partial_t\betae\\\no
&=\limsup_{\varepsilon\searrow0}-\theta_c\left(\int_Q|\partial_t\betae|^2-\frac 1 2\int_\Omega|\nabla\betae(t)|^2
+\frac 1 2\int_\Omega|\nabla\beta_0|^2-\int_Q\xie\partial_t\betae+\int_Q\partial_t\betae\right)\\\no
&=-\liminf_{\varepsilon\searrow0}\theta_c\left(\int_Q|\partial_t\betae^2|+\frac 12\int_\Omega|\nabla\betae(t)|^2
-\frac 12\int_\Omega|\nabla\beta_0|^2+\int_Q\xie\partial_t\betae-\int_Q\partial_t\betae\right)\\\no
&\leq -\theta_c\left(\int_Q\partial_t\beta^2-\frac 12\int_\Omega|\nabla\beta(t)|^2+\frac 12\int_\Omega|\nabla\beta_0|^2
-\int_Q\xi\partial_t\beta+\int_Q\partial_t\beta\right)\\\no
&=-\int_Q\theta\partial_t\beta,
\end{align}
by the fact that we have already identified the limit of \eqref{eqII} (written for $\varepsilon>0$, as $\varepsilon\searrow0$.
Note that we have exploited  weak lower semicontinuity of norms with \eqref{convIe}-\eqref{convIIe}, and the lower semicontinuity
of the indicator function.
Thus, we are in the position of proving the analogous of \eqref{sci1} and thus \eqref{sci2}-\eqref{sci3} easily follow,
concluding our proof.


\end{document}